\documentclass{article}
\usepackage[dvips]{graphicx}
\newcommand{\N}{\mathcal{N}}

\begin{document}
\title{$su_q(2)$-Invariant Schr\"odinger Equation of the
Three-Dimensional Harmonic Oscillator}
\author{M. IRAC-ASTAUD \\
Laboratoire de Physique Th\'eorique de la Mati\`ere
Condens\'ee,\\
Universit\'e
Paris VII, 2, Place Jussieu,\\ F-75251 Paris Cedex 05, France\\
e-mail : mici@ccr.jussieu.fr
\and
C.\ QUESNE \thanks{Directeur de recherches FNRS}\\
Physique Nucl\'eaire Th\'eorique et Physique
Math\'ematique,\\
Universit\'e Libre
de Bruxelles, Campus de la Plaine CP229,\\ Boulevard du Triomphe, B-1050 Brussels,
Belgium \\
e-mail : cquesne@ulb.ac.be}
\date{}
\maketitle
To appear in Lett.Math.Phys.
\begin{abstract}

We propose a $q$-deformation of the $su(2)$-invariant
Schr\"odinger equation of a spinless particle in a central potential,
which allows us not only to determine a
deformed spectrum and the corresponding eigenstates, as in other approaches,
but also to calculate the expectation values of some physically-relevant
operators. Here we consider the case of the isotropic harmonic oscillator
and of the quadrupole operator governing its interaction with an external
field. We obtain the  spectrum
and wave functions both for $q \in R^+$ and generic $q \in S^1$, and study the
effects of the $q$-value range and of the arbitrariness in the
$su_q(2)$ Casimir
operator choice. We then show that the quadrupole operator in $l=0$ states
provides
a good measure of the deformation influence on the wave functions and on the
Hilbert space spanned by them.
\end{abstract}



%
\section{Introduction}
Since the advent of quantum groups and quantum algebras~\cite{drinfeld, jimbo,
reshetikhin}, there has been a lot of interest in deformations of the harmonic
oscillator, since the latter plays a central role in the investigation of many
physical systems. Most studies were concerned with the one-dimensional
oscillator~\cite{arik, biedenharn,macfarlane}. Various $q$-deformed versions of
standard quantum mechanics in the Schr\"odinger representation were proposed for
the latter by using either the ordinary differentiation operator (see
e.g.~\cite{macfarlane}), or a $q$-differentiation one (see e.g.~\cite{kulish,
minahan,truong}).\par
%
%
The many-body problem for one-dimensional $q$-oscillators or, equivalently, the
deformation of the $N$-dimensional harmonic oscillator in cartesian coordinates
was also considered~\cite{floratos}, but relatively few works dealt with the
deformation of the same in radial and angular coordinates.
%
%
In one approach~\cite{jarvis}, only the radial problem was deformed via a purely
algebraic procedure based upon a $\mbox{so}_q(2,1) \oplus \mbox{so}(N)$
dynamical symmetry algebra. Only the $N=2$ case was dealt with in detail,
although extensions to higher dimension would in principle be feasible.\par
%
%
In another type of analysis, use was made of a differential calculus on the
$N$-dimensional non-commutative Euclidean space~\cite{carow91} to construct
and solve a $q$-deformed Schr\"odinger equation for the isotropic harmonic
oscillator~\cite{carow91,fiore,carow94,papp}. In such a setting, the underlying
SO($N$) symmetry gets replaced by an SO$_q$($N$) symmetry.
%
%
In still another study~\cite{chan}, the wave functions and energy spectrum
of the
$q$-isotropic oscillator were obtained in a rather indirect way from some
previous
results for the $q$-linear oscillator~\cite{finkelstein}. The latter were
derived by
replacing classical Poisson brackets by $q$-commutators (instead of standard
commutators).
%
%
In the last two approaches, the wave functions involve some non-commutative
objects: the variables of the quantum Euclidean space in the former, and the
elements of the SU$_q$(2) quantum group in the latter. This complicates the
calculation and interpretation of operator matrix elements.\par
%
%
The purpose of this Letter is to present an entirely different approach to the
$q$-deformation of the su(2)-invariant Schr\"odinger equation of a spinless
particle in a central potential, which allows us not only to determine a
deformed
spectrum as in other works, but also to easily calculate the expectation
values of
some physically-relevant operators, thereby evaluating the deformation influence
on the corresponding wave functions. Here we consider in detail the case of a
particle in an isotropic harmonic oscillator and of the quadrupole operator
governing its interaction with an external field.\par
%
%
In our approach, only the angular sector is deformed by using a
representation of the
su$_q$(2) quantum algebra on the two-dimensional sphere~\cite{rideau}. This
gives
rise to an appropriate change in the angular part of the scalar
product~\cite{irac98}, and to the substitution of the su$_q$(2) Casimir operator
eigenvalue for the su(2) one in the radial Schr\"odinger
equation~\cite{irac96}. The
latter step may be performed in various ways since there is no unique rule for
constructing the su$_q$(2) Casimir operator. Similarly, the deforming parameter
$q$ may be assumed either real and positive, or on the unit circle in the
complex
plane (but different from a root of unity), provided different scalar
products are
used~\cite{irac98}. We will study the effects of these two choices on the
solutions
of the radial Schr\"odinger equation.\par
%
%
In Section~\ref{sec:equation}, the su$_q$(2)-invariant Schr\"odinger equation of
the three-dimensional harmonic oscillator is introduced and solved. In
Section~\ref{sec:spectrum}, its spectrum is studied in detail for various
choices
of su$_q$(2) Casimir operators and $q$ ranges. The effect of the
deformation on the
corresponding wave functions is determined in Section~\ref{sec:quadrupole} by
calculating the quadrupole moment in $l=0$ states. Finally,
Section~\ref{sec:conclusion} contains the conclusion.\par
%
%
%
\section{ $su_q(2)$-Invariant Schr\"odinger Equation}
\label{sec:equation}
Let
\begin{equation}
  H_q = -\frac{\hbar^2}{2\mu} \left(\frac{\partial^2}{\partial r^2} +
\frac{2}{r}
  \frac{\partial}{\partial r} - \frac{C_q}{r^2}\right) + \frac{1}{2} \mu
\omega^2 r^2
  \label{eq:H}
\end{equation}
be the Hamiltonian of a $q$-deformed three-dimensional harmonic oscillator in
spherical coordinates $r$, $\theta$, $\phi$. Here $C_q$ is the $su_q(2)$ Casimir
operator, which we may take as
\begin{equation}
  C_q = J_+ J_- + \left[J_3 - \frac{1}{2}\right]_q^2 - \frac{1}{4},
\label{eq:C}
\end{equation}
where $[x]_q \equiv \left(q^x - q^{-x}\right)/\left(q - q^{-1}\right)$, and $q =
e^{w} \in R^+$ or $q = e^{{\rm i}w} \in S^1$ (but different from a
root of
unity). The operators $J_3$, $J_+$, $J_-$, satisfying the $su_q(2)$ commutation
relations
\begin{equation}
  \left[J_3, J_{\pm}\right] = \pm J_{\pm}, \qquad \left[J_+, J_-\right] =
  [2J_0]_q,
\end{equation}
are defined in terms of the angular variables by
\begin{eqnarray}
  J_3 & = & - {\rm i} \partial_{\phi}, \nonumber \\
  J_+ & = & - e^{{\rm i}\phi} \left( \tan(\theta/2) [T_1]_q\, q^{T_2}
          + \cot(\theta/2) q^{T_1}\, [T_2]_q\right), \nonumber \\
  J_- & = & e^{-{\rm i}\phi} \left( \cot(\theta/2) [T_1]_q\, q^{T_2}
          + \tan(\theta/2) q^{T_1}\, [T_2]_q\right),  \label{eq:J}
\end{eqnarray}
with $T_1 = - \frac{1}{2} \left(\sin \theta \partial_{\theta} - {\rm i}
\partial_{\phi}\right)$, $T_2 = - \frac{1}{2} \left(\sin \theta
\partial_{\theta} +
{\rm i} \partial_{\phi}\right)$~\cite{rideau,irac96}.\par
%
%
Instead of Equation~(\ref{eq:C}), we may alternatively use the operator
\begin{equation}
  C'_q = J_+ J_- + \left[J_3\right]_q \left[J_3-1\right]_q  \label{eq:C'}
\end{equation}
in Equation~(\ref{eq:H}), in which case the corresponding Hamiltonian will be
denoted by~$H'_q$.\par
%
%
The Hamiltonians $H_q$ and $H'_q$ remain invariant under $su_q(2)$ since they
commute with $J_3$, $J_+$, $J_-$, and they coincide with the Hamiltonian of the
standard three-dimensional isotropic oscillator when $q=1$. For
simplicity's sake,
we shall henceforth adopt units wherein $\hbar=\mu=~\omega~=~1$.\par
%
%
The $su_q(2)$-invariant Schr\"odinger equation
\begin{equation}
  H_q\, \psi_{nlmq}(r,\theta,\phi) = E_{nlq}\, \psi_{nlmq}(r,\theta,\phi)
\end{equation}
is separable and the corresponding wave functions can be written as
\begin{equation}
  \psi_{nlmq}(r,\theta,\phi) = {\cal R}_{nlq}(r) Y_{lmq}(\theta,\phi).
\label{eq:wf}
\end{equation}
\par
%
%
The $q$-spherical harmonics $Y_{lmq}(\theta,\phi)$, satisfying the equations
\begin{eqnarray}
  C_q Y_{lmq}(\theta,\phi) & = & C_q(l) Y_{lmq}(\theta,\phi), \qquad C_q(l) =
           \left[l+\frac{1}{2}\right]_q^2 - \frac{1}{4}, \\
  J_3 Y_{lmq}(\theta,\phi) & = & m Y_{lmq}(\theta,\phi),
\end{eqnarray}
are given by~\cite{rideau,irac98}
\begin{equation}\label{eq:R}
\begin{array}{ll}
 & Y_{lmq}(\theta,\phi)  =  {\cal N}_{lmq}\,
Q_{lq}\left(\cot^2(\theta/2)\right)
          R^l_{m0q}\left(\cot^2(\theta/2)\right) \cot^m(\theta/2)\, e^{{\rm
i}m\phi},
           \\
           &\\
 & {\cal N}_{lmq}  =  (-1)^l \left([2l+1]_q [l+m]_q!\right)^{1/2}\left(4\pi
          [l-m]_q!\right)^{-1/2},  \\
&Q_{lq}\left(\cot^2(\theta/2)\right)  = \prod_{k=0}^{l-1} \left(1 +
q^{2k-2l}
          \cot^2(\theta/2)\right)^{-1}, \\
& R^l_{m0q}\left(\cot^2(\theta/2)\right)  =  [l]_q!\, [l-m]_q! \sum_k
          \frac{\left(-\cot^2(\theta/2)\right)^k}{[k]_q!\, [l-m-k]_q!\,
[l-k]_q!\,
          [m+k]_q!},
\end{array}
\end{equation}
and $[x]_q! \equiv [x]_q [x-1]_q \ldots [1]_q$ if $x \in \N^+$, $[0]_q!
\equiv 1$, and
$\left([x]_q!\right)^{-1} \equiv 0$ if $x \in \N^-$. The functions
$Y_{lmq}(\theta,\phi)$, where $l=0$, 1, 2,~\ldots, and $m = -l$, $-l+1$,
\ldots,~$l$,
form an orthonormal set with respect to the scalar product~\cite{irac98}
\begin{equation}\label{eq:scalprod}
\begin{array}{ll}
  \langle l' m' | l m \rangle_q & \equiv  (q-q^{-1})(4\ln q)^{-1} \int_0^{\pi}
         d\theta\, \sin\theta \int_0^{2\pi} d\phi \nonumber \\
   & \mbox{} \times \Biggl(\overline{Y_{l'm'q^{\mp1}}(\theta,\phi)}\, \frac{1}
         {\sin^2(\theta/2) + q^{-2} \cos^2(\theta/2)}\,
q^{\sin\theta\partial_{\theta}
         - 1} Y_{lmq}(\theta,\phi) \nonumber \\
   & \mbox{} + \overline{Y_{l'm'q^{\pm1}}(\theta,\phi)}\,
\frac{1}{\sin^2(\theta/2)
         + q^2 \cos^2(\theta/2)}\, q^{- \sin\theta\partial_{\theta} + 1}
         Y_{lmq^{-1}}(\theta,\phi) \Biggr) \nonumber \\
  & =  \delta_{l',l} \delta_{m',m},
\end{array}
\end{equation}
where the upper (resp.\ lower) signs correspond to $q \in R^+$ (resp.\ $q
\in S^1$).
Note that the definitions~(\ref{eq:J}) and (\ref{eq:R})
slightly differ
from those given in~\cite{rideau}. With the present choice, in the $q\to1$ limit
they go over into the definitions of $su(2)$ generators and spherical
harmonics used
in most quantum mechanics textbooks (see e.g.~\cite{edmonds}).\par
%
%
The radial wave functions ${\cal R}_{nlq}(r) \equiv r^{-1} S_{nlq}(r)$ in
Equation~(\ref{eq:wf}) are the solutions of the radial equation
\begin{equation}
  \left(\frac{d^2}{dr^2} - \frac{C_q(l)}{r^2} - r^2 + 2 E_{nlq}\right)
S_{nlq}(r) = 0
\end{equation}
that satisfy the condition :
\begin{equation}
  S_{nlq}(0) = 0,
\end{equation}
and are square integrable with respect to the usual scalar product, i.e.,
\begin{equation}
  \langle n' l | n l \rangle_q \equiv \int_0^{\infty} dr\,
\overline{S_{n'lq}(r)}
  S_{nlq}(r) = \delta_{n',n}.
\end{equation}
They are given by
\begin{equation}
  S_{nlq}(r) = \sqrt{2 (n!)}\left( \Gamma\left(\alpha_l + n +
  \frac{1}{2}\right)\right)^{-1/2} e^{- \frac{1}{2} r^2} r^{\alpha_l}
L_n^{\alpha_l -
  \frac{1}{2}}(r^2), \label{eq:S}
\end{equation}
where $L_n^{\alpha_l - \frac{1}{2}}\left(r^2\right)$ is an associated Laguerre
polynomial~\cite{abramowitz}, and $\alpha_l$ is any one of the two solutions
\begin{equation}
  \alpha_{l\pm} = \frac{1}{2} \pm \lambda_q(l), \qquad \lambda_q(l) \equiv
  \sqrt{\frac{1}{4} + C_q(l)}, \label{eq:roots}
\end{equation}
of the equation
\begin{equation}
  \alpha_l (\alpha_l  - 1) = C_q(l), \label{eq:alpha-equation}
\end{equation}
provided it satisfies the condition $\alpha_l \in R^+$. The corresponding
energy
eigenvalues are
\begin{equation}
  E_{nlq} = 2n + \alpha_l + \frac{1}{2}, \qquad n = 0, 1, 2, \ldots.
\label{eq:E}
\end{equation}
The $\frac{1}{2}(N+1)(N+2)$-degeneracy of the isotropic oscillator energy
levels,
where $N = 2n+l$, is therefore lifted.\par
%
%
Similar results hold for the choice~(\ref{eq:C'}) for the Casimir operator,
the only
change being the substitution of $C'_q(l) = [l]_q [l+1]_q$ for $C_q(l)$. To
distinguish
the latter choice from the former, we shall denote all quantities referring
to it by
primed letters ($\lambda'_q(l)$, $\alpha'_l$, $E'_{nlq}$, ${\cal R}'_{nlq}(r)$,
$S'_{nlq}(r)$, \ldots).\par
%
%
\section{ Spectrum of the $su_q(2)$-Invariant Harmonic Oscillator}
\label{sec:spectrum}
In the present section, we will study the condition $\alpha_l \in R^+$ for the
existence of the radial wave functions~(\ref{eq:S}) and of the
corresponding energy
eigenvalues~(\ref{eq:E}), as well as the behaviour of the latter as
functions of $l$
and $q$ for the two choices (\ref{eq:C}), (\ref{eq:C'}) of Casimir
operators, and for
$q \in R^+$ or $q \in S^1$.\par
%
%
Let us first consider the case where $q = e^{w} \in R^+$. Since the
spectrum is
clearly invariant under the substitution $q \to q^{-1}$, we may assume
$q>1$, i.e.,
$w > 0$.\par
%
%
In the $C_q$ case, $\lambda_q^2(l) = \sinh^2((l+1/2)w) / \sinh^2w > 0$ for
$l=0$, 1, 2,~\ldots, hence both roots~(\ref{eq:roots}) of
Equation~(\ref{eq:alpha-equation}) are real and distinct. However, for $l\ne0$,
$\lambda_q(l) = \sinh((l+1/2)w) / \sinh w > 1$, showing that only
$\alpha_{l+}$ is positive and therefore admissible, whereas for $l=0$,
$\lambda_q(0) = [2 \cosh(w/2)]^{-1} < 1/2$ if $q\ne1$, so that $\alpha_{0+}$
and $\alpha_{0-}$ are both admissible. Note that in the undeformed case ($q=1$),
one gets $\lambda_1(0) = 1/2$, so that the root $\alpha_{0-}$ has then to be
discarded in accordance with known results. For $q\ne1$, the spectrum therefore
comprises the energy eigenvalues
\begin{eqnarray}
  E_{n0q\pm} & = & 2n + 1 \pm \frac{1}{2 \cosh(w/2)}, \label{eq:l=0-C-R} \\
  E_{nlq} & = & 2n + 1 + \frac{\sinh((l+1/2)w)}{\sinh w}, \qquad l =
1, 2, \ldots.
          \label{eq:l-C-R}
\end{eqnarray}
The appearance of an additional $l=0$ level was already observed in the
deformed Coulomb potential case~\cite{irac96}.\par
%
%
In the $C'_q$ case, $\lambda_q^{\prime2}(l) = [l]_q [l+1]_q + \frac{1}{4} \ge
\frac{1}{4}$ for $l=0$, 1, 2,~\ldots, so that both roots $\alpha'_{l+}$ and
$\alpha'_{l-}$ are real and distinct, but only the former is positive, hence
admissible. The spectrum therefore comprises the same levels as in the
undeformed case, their energies being now
\begin{equation}
  E'_{nlq} = 2n + 1 + \frac{[4 \sinh(lw) \sinh((l+1)w) +
\sinh^2w]^{1/2}}
  {2\sinh w}, \quad l = 0, 1, 2, \ldots  \label{eq:l-C'-R}
\end{equation}
Note that the energy of the $l=0$ states is left undeformed :
\begin{equation}
  E'_{n0q} = E_{n0} = 2n + \frac{3}{2}.  \label{eq:l=0-C'}
\end{equation}
\par
%
%
Expanding the right-hand sides of Eqs.~(\ref{eq:l=0-C-R}), (\ref{eq:l-C-R}),
and~(\ref{eq:l-C'-R}) into powers of $w$ shows that in the neighbourhood of
$q=1$, i.e., $w=0$,
\begin{equation}\label{eq:E-expansion}
\begin{array}{lll}
 & E_{n0q\pm} & \simeq  2n + 1 \pm \frac{1}{2} \mp \frac{1}{16} w^2 \left(1 -
           \frac{5}{48} w^2 + \cdots\right), \\
           &&\\
  &E_{nlq} & \simeq  2n + l + \frac{3}{2} + \frac{1}{48} (2l-1) (2l+1)
(2l+3) w^2 \\
           & & \times \left\{1 +\frac{1}{240} [12l(l+1)-25] w^2 + \cdots\right\},
  \quad l = 1, 2, \ldots,  \\
  &&\\
  &E'_{nlq} & \simeq  2n + l + \frac{3}{2} + \frac{l(l+1)}{6(2l+1)} w^2
\biggl\{
           2l(l+1) - 1  \\
  &&  \mbox{} + \frac{24l^3(l+1)^3 - 56l^2(l+1)^2 - 10l(l+1) + 7}{60
[4l(l+1)+1]}
           w^2 + \cdots\biggr\},
           \quad  l = 0, 1, 2, \ldots.
\end{array}
\end{equation}
Hence, for $l\ne0$, $E_{nlq}$ and $E'_{nlq}$ are increasing functions of
$w$ in
the neighbourhood of $w=0$, whereas for $l=0$, $E_{n0q+}$ and $E_{n0q-}$ have
opposite behaviours, while $E'_{n0q}$ is independent of $w$. Moreover, for a
given $w$ value, the influence of the deformation increases with $l$. Such
trends are confirmed by Figure~1, where the first few lowest eigenvalues
$E_{0lq}$ are plotted in terms of $w$. For the $l$ and $w$ values
considered,
$E'_{0lq}$ cannot be distinguished from $E_{0lq}$ for $l\ne0$, or $E_{00q+}$ for
$l=0$.
\begin{figure}
\center{ \includegraphics[angle=270,scale=.4]{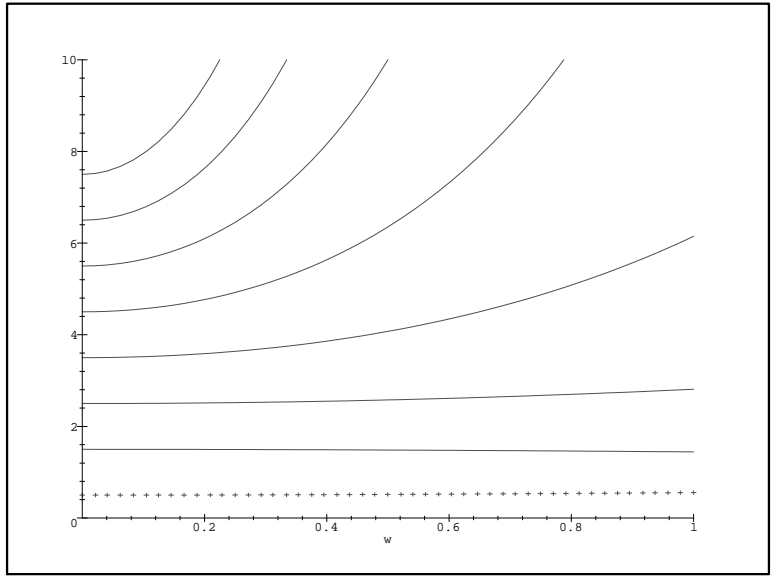}}
  \caption{Spectrum of $n=0$ states in terms of $w$ for $q =
e^{w} \in {\mbox R}^+$. The solid lines correspond to $E_{00q+}$ and $E_{0lq}$, $l=1$, 2,
\ldots,~6,
and the  crossed line to $E_{00q-}$.}
\end{figure}
%
%

Let us next consider the case where $q = e^{{\rm i}w} \in S^1$. Owing to the
invariance of the spectrum under the substitution $q \to q^{-1}$, we may now
assume $0 < w < \pi$.\par
%
%
{}For $C_q$, $\lambda_q^2(l) = \sin^2((l+1/2)w)/\sin^2w$ vanishes for
$w =
k\pi/(l+1/2)$, $k=1$, 2, \ldots,~$l$, but these $w$ values are in any case
excluded as roots of unity. Hence, for the $w$ values considered here,
Equation~(\ref{eq:alpha-equation}) has two real, distinct roots. If
$\lambda_q(l) =
|\sin((l+1/2)w)/\sin w| < 1/2$, both roots $\alpha_{l+}$ and
$\alpha_{l-}$ are admissible, whereas if $\lambda_q(l) \ge 1/2$, only
$\alpha_{l+}$ is so. The conditions can be reformulated in terms of
\begin{equation}
  \gamma_q(l) = 4 \sin^2w C_q(l) = \frac{1}{2} [- 4 \cos((2l+1)w) +
\cos2w
  + 3].
\end{equation}
One finds two admissible roots if $\gamma_q(l) < 0$, but only one root
$\alpha_{l+}$ if $\gamma_q(l) \ge 0$. For instance, for $l=0$, there is a single
eigenvalue for any $w$ value,
\begin{equation}
  E_{n0q} = 2n + 1 + \frac{1}{2\cos(w/2)},  \label{eq:l=0-C-S}
\end{equation}
while for $l=1$, there are either two or one eigenvalues, which are given by
\begin{equation}
\begin{array}{llll}
  &E_{n1q\pm} &= 2n + 1 \pm \frac{4\cos^2(w/2) - 1}{2\cos(w/2)},   &
          \mbox{if\ } \frac{-7-\sqrt{17}}{16} < \cos w <
\frac{-7+\sqrt{17}}{16}, \\
&&&\\
  &E_{n1q} &= 2n + 1 + \frac{4\cos^2(w/2) - 1}{2\cos(w/2)},   &
          \mbox{if\ } \cos w \le \frac{-7-\sqrt{17}}{16}  \\
         & &&\\
  && & \mbox{or\ } \cos w \ge \frac{-7+\sqrt{17}}{16},
\end{array}
\end{equation}
respectively.\par
%
%
By proceeding in a similar way for $C'_q$, one finds in terms of
\begin{equation}
  \gamma'_q(l) = 4 \sin^2w C'_q(l) = 4 \sin((l+1)w) \sin(lw)
\end{equation}
that there are two admissible roots $\alpha'_{l+}$, $\alpha'_{l-}$ if $-
\sin^2w <
\gamma'_q(l) < 0$, only one $\alpha'_{l+}$ if either $\gamma'_q(l) \ge 0$ or
$\gamma'_q(l) = - \sin^2w$ (in which case $\alpha'_{l+} = 1/2$), or none if
$\gamma'_q(l) < - \sin^2w$. For instance, for $l=0$, there is a single
eigenvalue (\ref{eq:l=0-C'}) for any $w$ value and it coincides with the
undeformed one  , while for $l=1$, there are two or one eigenvalues if
$\cos w \ge
- 1/8$,
\begin{equation}
\begin{array}{llll}
  &E'_{n1q\pm}  &=  2n + 1 \pm \frac{1}{2} \sqrt{1 + 8\cos w}, &
\mbox{if\ }
           - \frac{1}{8} < \cos w < 0, \\
           &&\\
  &E'_{n1q} &=  2n + 1 + \frac{1}{2} \sqrt{1 + 8\cos w}, &\mbox{if\ }
           \cos w = - \frac{1}{8} \mbox{\ or\ } \cos w \ge 0,
\end{array}
\end{equation}
and none if $\cos w < -1/8$.\par
%
%
In both the $C_q$ and $C'_q$ cases, similar results can be derived for
higher $l$
values, and close enough to $w=0$, for given $n$ and $l$ values one
always finds
a single eigenvalue going into the undeformed one, $E_{nl} = 2n + l + 3/2$,
for $q
\to 1$. For all $l$ values, the expansion of this eigenvalue into powers of
$w$
can be obtained from Equation~(\ref{eq:E-expansion})
by substituting ${\rm i} w$ for $w$. Hence, in a small enough
neighbourhood of
$w=0$, for $l\ne0$, $E_{nlq}$ and $E'_{nlq}$ are decreasing functions of
$w$ ,
whereas $E_{n0q}$ is increasing and $E'_{n0q}$ remains constant.\par
%
%
\begin{figure}
\center{ \includegraphics[angle=270,scale=.4]{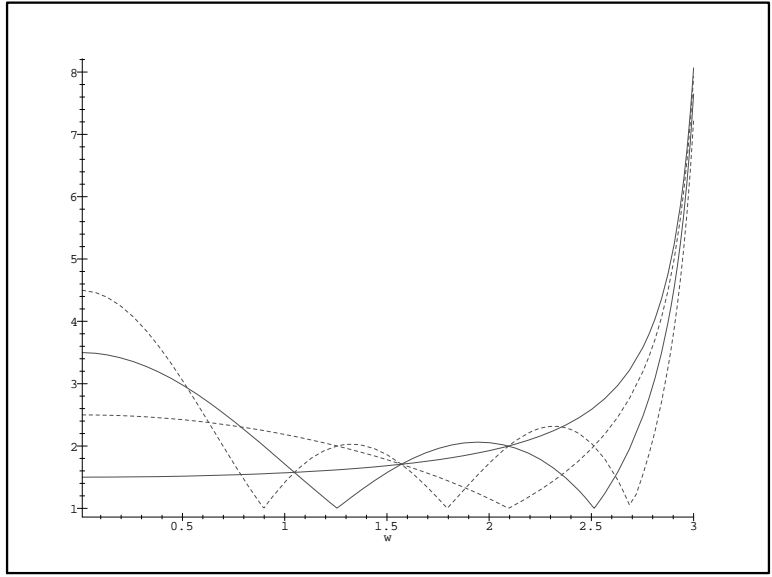}}
  \caption{ Spectrum of $n=0$ states in terms of $w$ for $q =
e^{{\rm i}w}
\in S^1$. The solid and dashed lines correspond to $E_{0lq}$ (or
$E_{0lq+}$ if there exists another eigenvalue with $n=0$ and the same $l$),
where $l=0, 2$,  and $l=1, 3$,  respectively.}
\end{figure}
The first few eigenvalues $E_{0lq}$, going into $E_{0l}$ when $q \to 1$, are
displayed on Figure~2 for $0 < w < \pi$. One should remember that for
some $l$
and $w$ values, there may exist other eigenvalues, which are not plotted
on the
figure, and that the discrete set of points, where $w$ is equal to a
root of unity,
is excluded. Such is the case, in particular, of the points where $E_{0lq}$
($l\ne0$)
takes its minimal value~1, and of the point $w = \pi/2$, where $E_{0lq} = 1 +
(1/\sqrt{2})$ for any $l$ value. The influence of the deformation on the
spectrum is
rather striking. For high $w$ values, the levels get mixed in a very
complicated
way. It is remarkable that in the neighbourhood of $w = \pi/2$, one obtains a
spectrum very close to one with equidistant, infinitely-degenerate levels.\par
%
%
In the $C'_q$ case, the situation is still more complex as some $l$ values may
disappear on some intervals. In the neighbourhood of $w=0$, however, as
in the
real $q$ case, $E'_{nlq}$ cannot be distinguished from $E_{nlq}$.\par
%
%
\section{ Quadrupole Moment in $l=0$ States}
\label{sec:quadrupole}
The purpose of the present section is to study the effect of the
deformation on the
wave functions of the $su_q(2)$-invariant harmonic oscillator, given in
Equation~(\ref{eq:wf}), (\ref{eq:R}), and~(\ref{eq:S}), by
determining
the variation with $w$ of the expectation value of some physically-relevant
operator. For the latter, we choose the electric quadrupole moment
operator, and we
consider the quadrupole moment in a state with definite $n$, $l$ values,
which is
defined conventionally as~\cite{edmonds}
\begin{equation}
  Q_{nlq} = \langle n l l| (3z^2 - r^2) | n l l\rangle_q.  \label{eq:Q}
\end{equation}
Equation~(\ref{eq:Q}) corresponds to the choice $C_q$ for the $su_q(2)$ Casimir
operator. When using instead $C'_q$, the quadrupole moment will be denoted by
$Q'_{nlq}$.\par
%
%
The undeformed counterpart of $Q_{nlq}$ and $Q'_{nlq}$ is given by
\begin{equation}
  Q_{nl} = \langle n l| r^2 |n l\rangle \langle l l| \left(3 \cos^2\theta - 1\right))| l
  l\rangle = \left(2n + l + \frac{3}{2}\right) \left(- \frac{2l}{2l+3}\right).
  \label{eq:Q-undef}
\end{equation}
It vanishes for $l=0$ as a result of the familiar selection rule for the angular
momenta $l$ and 2 coupling. Nonvanishing values of $Q_{n0q}$ or $Q'_{n0q}$ will
therefore be a direct measure of the effect of the deformation. Note that
since in
the $C_q$ case, there are two energy eigenvalues with $l=0$ and a given $n$
value
for $q \in R^+$ (see Equation~(\ref{eq:l=0-C-R})), we have to distinguish the
corresponding quadrupole moments by a $\pm$ subscript.\par
%
%
As in Equation~(\ref{eq:Q-undef}), $Q_{n0q}$ can be factorized into radial and
angular matrix elements,
\begin{equation}
  Q_{n0q} = \langle n 0| r^2 |n 0\rangle_q \langle 00| \left(3 \cos^2\theta
- 1\right))|
  00\rangle_q.
\end{equation}
For any $l$ value, the former is simply obtained by replacing $l$ by
$\alpha_l - 1$
in the undeformed radial matrix element. Hence
\begin{equation}
  \langle n 0| r^2 |n 0\rangle_q = 2n + \alpha_0 + \frac{1}{2}.
\label{eq:def-radial}
\end{equation}
The calculation of the latter is more complicated as it implies the use of the
deformed angular scalar product~(\ref{eq:scalprod}),
\begin{equation}
  \langle 00| (3 \cos^2\theta - 1)| 00\rangle_q = \frac{3 (q-q^{-1})}
  {16 \pi \ln q} \left({\cal I}_q + {\cal I}_{q^{-1}}\right) - 1,
  \label{eq:def-angular}
\end{equation}
where
\begin{equation}
  {\cal I}_q = \int_{0}^{\pi} d\theta\, \sin\theta \int_{0}^{2\pi} d\phi\,
  \frac{1}{\sin^2(\theta/2) + q^{-2} \cos^2(\theta/2)}\, q^{\sin\theta
  \partial_{\theta} - 1} \cos^2\theta.
\end{equation}
This integral can be easily performed by making the changes of variables $z
= \rho
e^{{\rm i}\phi}$, $\overline{z} = \rho e^{-{\rm i}\phi}$, $\rho =
\cot(\theta/2)$,
and $\eta = \rho^2$~\cite{rideau, irac98}. One gets
\begin{equation}
  {\cal I}_q = \frac{4 \pi}{ q} \int_{0}^{\infty} d\eta\,
\frac{(1-q^{-2}\eta)^2}{(1+\eta)
  (1+ q^{-2}\eta)^3} = 8 \pi \frac{q+q^{-1}}{(q-q^{-1})^2}
  \left(\frac{q+q^{-1}}{q-q^{-1}} \ln q - 1\right).  \label{eq:I}
\end{equation}
\par
%
%
Introducing Equation~(\ref{eq:I}) into Equation~(\ref{eq:def-angular}), we
finally
obtain
\begin{equation}
  \langle 00| (3 \cos^2\theta - 1)| 00\rangle_q = \frac{2 (q^2+4+q^{-2})}
  {(q-q^{-1})^2} - \frac{3 (q+q^{-1})}{(q-q^{-1}) \ln q}.
\end{equation}
Hence
\begin{equation}\label{eq:def-angular-bis}
\begin{array}{lll}
  \langle 00| (3 \cos^2\theta - 1)| 00\rangle_q & =  \frac{2
         \cosh^2w + 1}{\sinh^2w} - \frac{3 \cosh w}{w
\sinh w}, \quad
         &\mbox{if\ } q = e^{w} \in R^+, \nonumber \\
         &&\\
  & =   - \frac{2 \cos^2w + 1}{\sin^2w} + \frac{3 \cos w}{w
\sin w},
        \quad  &\mbox{if\ } q = e^{{\rm i}w} \in S^1.
\end{array}
\end{equation}
The deformed quadrupole moment $Q_{n0q}$ is therefore an even function of
$w$, so that we may again restrict ourselves to $0 < w < \infty$ or
$0 < w <
\pi$ according to whether $q$ is real or complex.\par
%
%
The deformed radial matrix element~(\ref{eq:def-radial}) (or its counterpart for
$Q'_{n0q}$), being just equal to the corresponding energy eigenvalue in the
units
used, varies in the same way with $w$. From Equations~(\ref{eq:l=0-C-R}),
(\ref{eq:l=0-C'}) and~(\ref{eq:l=0-C-S}), it follows that for increasing
$w$, it
decreases (resp.\ increases) from $2n + \frac{3}{2}$ (resp.\ $2n +
\frac{1}{2}$) to
$2n+1$ for $C_q$, $q \in R^+$, and $\alpha_{0+}$ (resp.\ $\alpha_{0-}$),
increases
from $2n + \frac{3}{2}$ to $+\infty$ for $C_q$, $q \in S^1$, and remains
constant
for $C'_q$, $q \in R^+$ or $q \in S^1$. The effect of the deformation is
therefore
not significant, except in the case of $C_q$ when $q \in S^1$.\par
%
%
{}For real $q$, the deformed angular matrix element, given in
Equation~(\ref{eq:def-angular-bis}), increases from 0 to 2 when $w$ goes from
0 to $+\infty$. Hence, it is obvious that $Q_{n0q-}$ and $Q'_{n0q}$ are
increasing
positive functions of $w$. It can be checked that the same is true for
$Q_{n0q+}$, the variation of the angular part of the matrix element
dominating that
of the radial one.
 On Figure~3, $Q_{n0q+}$, $Q_{n0q-}$, and $Q'_{n0q}$ are
displayed
in terms of $w$ for $n=0$ and $n=1$.
%
%
 One can see some effect of the choice of
Casimir operator and $\alpha_0$ root, but it becomes significant only for a very
large deformation. It should be stressed that the undeformed quadrupole moments
in $l\ne0$ states and the deformed ones in $l=0$ states have opposite
signs. Both
become comparable in absolute value for $w \ge 1$.
\begin{figure}
\center{ \includegraphics[angle=270,scale=.4]{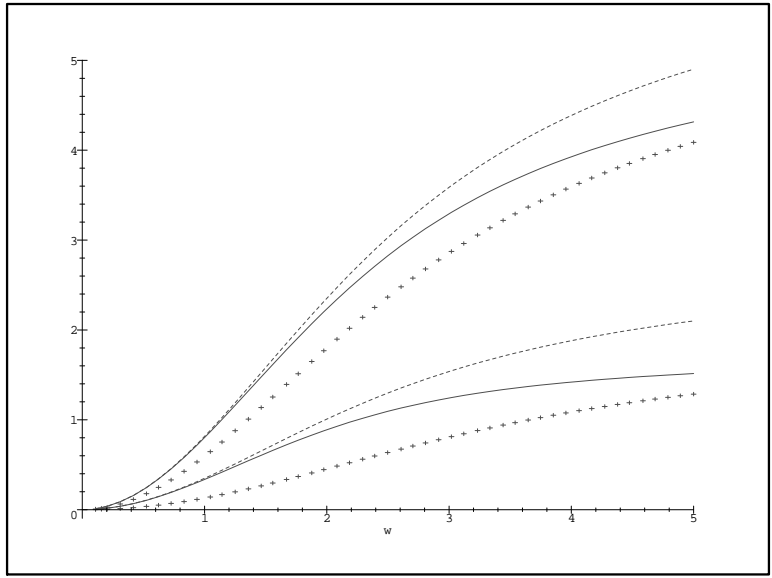}}
  \caption{  Quadrupole moment of $l=0$ states in terms of $w$ for $q =
e^{w} \in {\mbox R}^+$. The solid, crossed, and dashed lines correspond to $Q_{n0q+}$,
$Q_{n0q-}$, and $Q'_{n0q}$, respectively, for $n=0$ (three lowest curves)
and $n=1$
(three highest ones).}
\end{figure}

{}For complex $q$, the deformed angular matrix
element~(\ref{eq:def-angular-bis})
decreases from 0 to $-\infty$ when $w$ goes from 0 to $\pi$. Hence, both
$Q_{n0q}$ and $Q'_{n0q}$ are decreasing negative functions of $w$. On
Figure~4,
they are displayed in terms of $w$ for $n=0$ and $n=1$. Apart from the sign,
which is now the same as that of $Q_{nl}$ for $l\ne0$, the conclusions remain
similar to those for the real $q$ case.
\begin{figure}
\center{ \includegraphics[angle=270,scale=.4]{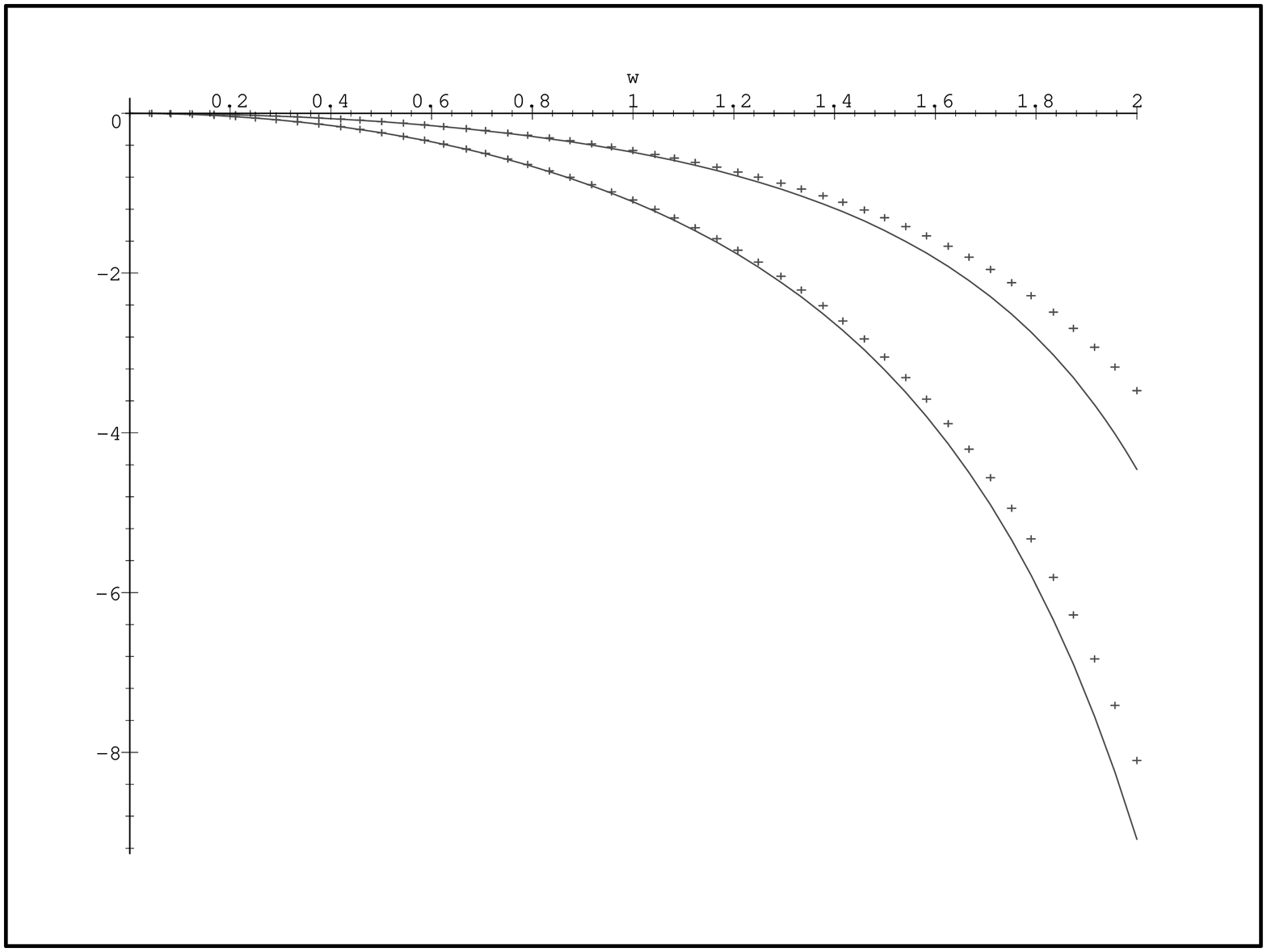}}
  \caption{   Quadrupole moment of $l=0$ states in terms of $w$ for $q =
e^{{\rm i}w} \in S^1$. The solid and crossed lines correspond to $Q_{n0q}$ and
$Q'_{n0q}$, respectively, for $n=0$ (two highest curves) and $n=1$ (two lowest
ones).}
\end{figure}
%
%
\section{Conclusion}
\label{sec:conclusion}

In the present Letter, we did show that as those of the free particle and of the
Coulomb potential~\cite{irac96}, the $su_q(2)$-invariant
Schr\"odinger equation
of the three-dimensional harmonic oscillator can be easily solved not only
for $q \in R^+$, but also for generic $q \in S^1$. It is worth stressing
that we have been
working in the framework of the usual Schr\"odinger equation (i.e., with no
non-commuting objects contrary to some other
approaches~\cite{carow91,fiore,carow94,papp,chan}), but with wave functions
belonging to a Hilbert space different from the usual one, since the
angular part of
the scalar product has been modified when going from $su(2)$ to
$su_q(2)$~\cite{irac98}.\par
%
%
In the real $q$ case, we did show that the spectrum is rather similar to the
undeformed one, except that the energy levels are no more equidistant and that
their degeneracy is lifted. For a given $n$ value, the spacing between adjacent
levels corresponding to $l$ and $l+1$, respectively, increases with $l$ and with
the deformation. In addition, there appears a supplementary series of $l=0$
levels
when the Casimir operator $C_q$ is used. Apart from this, for small
deformations,
the results are rather insensitive to the choice made for the Casimir
operator.\par
%
%
In the complex $q$ case, we did show that the spectrum is more complicated
as for
$l\ne0$ and any $n$ value, there may exist 0, 1, or 2 levels according to the
deformation. The existence or inexistence of levels is also rather
sensitive to the
choice made for the Casimir operator. Close enough to $q=1$, there however
always
exists a single level going into the undeformed one for $q \to 1$. In that
region, the
spacing between adjacent levels corresponding to $l$ and $l+1$,
respectively, now
decreases with $l$ and with the deformation.\par
%
%
The closeness of our approach to the standard one did also allow us to study the
effect of the deformation on the wave functions and the Hilbert space spanned by
them. We did establish that it is rather strong as the quadrupole moment in the
$l=0$ states, which vanishes in the undeformed case, now assumes a positive
(resp.\ negative) value for $q \in R^+$ (resp.\ $q \in S^1$) irrespective
of the
Casimir operator used.\par
%
%

%
%
\end{document}